\documentclass[11pt]{article}

\usepackage{cmap}
\usepackage[T1]{fontenc}
\usepackage[utf8]{inputenc}
\usepackage[english]{babel}
\usepackage{amsmath,amssymb,amsfonts,amsthm}
\usepackage[numbers,sort&compress]{natbib}
\usepackage{xcolor}
\usepackage{url}
\usepackage{hyperref}

\DeclareTextSymbolDefault{\textsection}{OMS}

\theoremstyle{plain}
\newtheorem{theorem}{Theorem}
\newtheorem*{lemma}{Lemma}
\theoremstyle{remark}
\newtheorem{remark}{Remark}
\newtheorem*{example}{Example}

\hypersetup{
  colorlinks=true,
  allcolors=black,
  pdftitle={The Fold Singularity of Solutions to the System Equivalent to the Nonlinear Stationary Transonic Equation},
  pdfauthor={Azamat Mavletovich Shavlukov}
}

\title{The Fold Singularity of Solutions to the System Equivalent to the Nonlinear Stationary Transonic Equation}
\author{Azamat Mavletovich Shavlukov\\[3pt]
\parbox{0.95\textwidth}{\centering\small Institute of Mathematics with Computing Centre, Ufa Federal Research Centre of the Russian Academy of Sciences}\\
\small 112 Chernyshevsky Street, Ufa 450008, Russia\\
\small \href{mailto:aza3727@yandex.ru}{aza3727@yandex.ru}\\
\small ORCID: \href{https://orcid.org/0000-0003-1779-839X}{0000-0003-1779-839X}}
\date{}

\begin{document}
\maketitle

\begin{abstract}
We study the fold singularity of solutions to the quasilinear system equivalent to the nonlinear stationary transonic equation. The hodograph transformation and the introduction of a potential reduce this system to the linear Tricomi equation. We consider analytic solutions and focus on points of gradient catastrophe whose images in the hodograph plane lie on the line where the Tricomi equation changes type. We use invertible analytic changes of variables to express solutions of the quasilinear system implicitly in terms of the roots of the canonical fold equation. We then establish the genericity of the fold on a compact neighborhood in the hodograph plane: solutions satisfying the specified transversality condition form an open dense subset of the space of analytic solutions of the Tricomi equation restricted to this compact set. The topology is that of uniform convergence of all derivatives on the compact set. Density is proved using a family of analytic solutions depending on two real parameters and Sard's theorem. We also show that density fails in the strong Whitney topology on the space of solutions defined on an open rectangle.
\end{abstract}

\noindent\textbf{Keywords:} gradient catastrophe, singularity theory of differentiable mappings, fold, transonic flow, Tricomi equation.

\medskip
\noindent\textbf{Mathematics Subject Classification (2020):} Primary 35M10; Secondary 58K05, 76H05.

\bigskip

\section*{Introduction}

Applications of singularity theory of differentiable mappings to solutions of quasilinear systems go back to work carried out in the 1990s by A.~Kh.~Rakhimov, a member of V.~I.~Arnold's school \cite{rakh92,rakh93}. Arnold was one of the founders of the mathematical theory of catastrophes. Subsequent studies have investigated generic singularities of solutions to various first-order quasilinear systems, including equations of gas dynamics, hydrodynamics, and nonlinear geometrical optics \cite{dub06,dgk09,ko13,ko17,ko20,bogaevsky,kud01,gar10,ks1995,suleimanov2009,sul12,chap24,ss2020,ss22,ss24,mss24,shav,mss2026,sa}.

In this paper we study the fold singularity of solutions to the quasilinear system equivalent to the nonlinear stationary transonic equation. The hodograph transformation and the subsequent introduction of a potential reduce this system to the linear Tricomi equation. We consider analytic solutions of the Tricomi equation and focus on a gradient-catastrophe point with finite coordinates on the sonic line, where the nonlinear stationary transonic equation changes type. Near this point, solutions of the quasilinear system are expressed in terms of solutions of the canonical fold equation. The genericity of the fold is established on a compact neighborhood in the hodograph plane.

\section{Description of the fold singularity}\label{sec-description}

Consider the nonlinear stationary transonic equation \cite[Eq.~(65)]{lin1948}
\begin{equation}
 w_xw_{xx}+w_{yy}=0.
 \label{te}
\end{equation}
In gas dynamics, $w=w(x,y)$ is the normalized potential for the velocity perturbation in the $(x,y)$-plane
\cite[\S\S\,2, 3, formulas~(4), (9)]{lin1948}.

Equation \eqref{te} is elliptic when $w_x>0$ and hyperbolic when $w_x<0$; it changes type across $w_x=0$.

Set
\[
 u=w_x,\qquad v=w_y
\]
to rewrite \eqref{te} as the equivalent quasilinear system
\begin{equation}
 \begin{cases}
 uu_x+v_y=0,\\
 u_y=v_x.
 \end{cases}
 \label{tls}
\end{equation}
We now follow the approach of \cite{ks1995,ss2020,mss2026}.
Using the derivative formulas for the hodograph transformation
(see, e.g., \cite[Chap.~1, \S\,4]{ry}),
\begin{equation}
 \begin{gathered}
 u_x=Jy_v,\quad u_y=-Jx_v,\quad v_x=-Jy_u,\quad v_y=Jx_u,\\
 J=u_xv_y-u_yv_x,\qquad j=x_uy_v-x_vy_u,\qquad J=j^{-1}
 \end{gathered}
 \label{g}
\end{equation}
we transform the quasilinear system \eqref{tls} into the linear system
\begin{equation}
 \begin{cases}
 uy_v+x_u=0,\\
 x_v=y_u.
 \end{cases}
 \label{hodlinear}
\end{equation}
Introduce a potential $B=B(u,v)$ by
\begin{equation}
 x=B_u,\qquad y=B_v.
 \label{po1}
\end{equation}
The linear system \eqref{hodlinear} then reduces (see, e.g., \cite[Chap.~1, \S\,4]{bers2016}) to the linear Tricomi equation of mixed type
\begin{equation}
 B_{uu}+uB_{vv}=0.
 \label{Tre}
\end{equation}
Equation \eqref{Tre} is elliptic for $u>0$ and hyperbolic for $u<0$; its type changes across $u=0$.

The connection between the Tricomi equation \eqref{Tre} and transonic flow was apparently first established in \cite{frankl1944,frankl1945}.

By \eqref{po1} and \eqref{Tre}, the Jacobian of the hodograph transformation \eqref{g} is
\begin{equation}
 j=-uB_{vv}^2-B_{uv}^2.
 \label{j}
\end{equation}

Suppose that $j(u_*,v_*)=0$, and set
\[
 x_*=B_u(u_*,v_*),\qquad y_*=B_v(u_*,v_*).
\]
We assume that $(x_u^2+x_v^2+y_u^2+y_v^2)|_{(u_*,v_*)}\ne0$. By the derivative formulas \eqref{g}, as $(x,y)$ approaches the finite gradient-catastrophe point $(x_*,y_*)$, at least some first derivatives of $u$ and $v$ become unbounded in absolute value. The functions $u$ and $v$ themselves tend to the finite values $u_*$ and $v_*$.

Consider the case $u_*=0$.
For any functions $f_0$ and $f_1$ analytic near $v_*$, the Cauchy data \cite[\S\,1, formula~(2)]{korobeinik1962}
\[
 B(0,v)=f_0(v),\qquad B_u(0,v)=f_1(v)
\]
locally determine an analytic solution $B(u,v)$ of \eqref{Tre} near $(0,v_*)$.
Set $\Delta v=v-v_*$ and expand the Cauchy data in convergent series near $v_*$:
\[
 f_0(v)=\sum_{j\ge0}b_{0j}\Delta v^j,\qquad
 f_1(v)=\sum_{j\ge0}b_{1j}\Delta v^j,
\]
\[
 b_{0j}=\frac{f_0^{(j)}(v_*)}{j!},\qquad
 b_{1j}=\frac{f_1^{(j)}(v_*)}{j!}.
\]
The resulting solution of \eqref{Tre} has a convergent expansion near $(0,v_*)$:
\[
 B(u,v)=\sum_{i,j\ge0}b_{ij}u^i\Delta v^j,\qquad
 b_{ij}=\left.\frac{1}{i!j!}
 \frac{\partial^{i+j}B}{\partial u^i\partial v^j}\right|_{(0,v_*)}.
\]
Substitution into \eqref{Tre} yields the recurrence relations
\begin{equation}
 b_{2,j}=0,\qquad
 b_{i+3,j}=-\frac{(j+2)(j+1)}{(i+3)(i+2)}b_{i,j+2},
 \qquad i,j\ge0.
 \label{k}
\end{equation}
These recurrence relations determine all remaining coefficients from the two sequences $b_{0j}$ and $b_{1j}$.
In particular,
\[
 b_{30}=-\frac{b_{02}}3,\qquad
 b_{31}=-b_{03},\qquad
 b_{40}=-\frac{b_{12}}6.
\]

The Cauchy--Kovalevskaya theorem \cite[Chap.~I, \S\,7, Parts~3, 4]{courant1964} provides an alternative proof of local existence for \eqref{Tre} near $(0,v_*)$ when the following analytic Cauchy data are prescribed:
\[
\begin{aligned}
	B(0,v)&=b_{00}+y_*\Delta v+b_{02}\Delta v^2
	+\Delta v^3h_0(v),\\
	B_u(0,v)&=x_*+b_{11} \Delta v+b_{12}\Delta v^2+\Delta v^3h_1(v),
\end{aligned}
\]
where $h_0$ and $h_1$ are arbitrary functions analytic near $v=v_*$. Smooth, nonanalytic Cauchy data do not in general ensure that a solution extends smoothly from the hyperbolic region $u<0$ across $u=0$ into the elliptic region $u>0$ \cite[Theorem~4.1]{matsuzawa1972}.

The condition $j(0,v_*)=0$ implies $b_{11}=0$. Equation \eqref{Tre} also gives $b_{20}=0$, whereas the conditions obtained so far do not imply $b_{02}=0$. We therefore consider the case
\[
 b_{02}=\frac12f_0''(v_*)\ne0.
\]
Section~2 establishes genericity in the sense defined there.

Since $j_u(0,v_*)=-B_{vv}(0,v_*)^2\ne0$, the equation $j=0$ defines an analytic curve $u=\gamma(v)$ near $(0,v_*)$. Along this curve, \eqref{j} gives $u=-(B_{uv}/B_{vv})^2\le0$. The tangent line to the curve at $(0,v_*)$ is therefore $u=0$.

By \eqref{j}--\eqref{k}, the solution of \eqref{Tre} has the expansion
\begin{equation}
 \begin{aligned}
 B={}&b_{00}+x_*u+y_*\Delta v+b_{02}\Delta v^2
 -\frac{b_{02}}3u^3+b_{12}u\Delta v^2+b_{03}\Delta v^3\\
 &-\frac{b_{12}}6u^4-b_{03}u^3\Delta v
 +b_{13}u\Delta v^3+b_{04}\Delta v^4
 +\sum_{i+j\ge5}b_{ij}u^i\Delta v^j.
 \end{aligned}
 \label{bsol}
\end{equation}
Set $\Delta x=x-x_*$ and $\Delta y=y-y_*$. By \eqref{po1} and \eqref{bsol},
\begin{align}
 \Delta x={}&-b_{02}u^2+b_{12}\Delta v^2
 -\frac23b_{12}u^3-3b_{03}u^2\Delta v+b_{13}\Delta v^3
 +\sum_{i+j\ge4}x_{ij}u^i\Delta v^j,\notag\\
 \Delta y={}&2b_{02}\Delta v+2b_{12}u\Delta v+3b_{03}\Delta v^2\notag\\
 &-b_{03}u^3+3b_{13}u\Delta v^2+4b_{04}\Delta v^3
 +\sum_{i+j\ge4}y_{ij}u^i\Delta v^j.
 \label{po2}
\end{align}

Since $b_{02}\ne0$, the second equation in \eqref{po2} can be solved for $\Delta v$ by the implicit function theorem. The resulting analytic function $\Delta v=\Delta v(u,\Delta y)$ has the expansion
\begin{equation}
 \begin{aligned}
 \Delta v={}&\frac{\Delta y}{2b_{02}}
 -\frac{3b_{03}}{8b_{02}^3}\Delta y^2
 +\left(\frac{9b_{03}^2}{16b_{02}^5}
 -\frac{b_{04}}{4b_{02}^4}\right)\Delta y^3
 -\frac{b_{12}}{2b_{02}^2}u\Delta y\\
 &+\left(\frac{9b_{12}b_{03}}{8b_{02}^4}
 -\frac{3b_{13}}{8b_{02}^3}\right)u\Delta y^2
 +\frac{b_{12}^2}{2b_{02}^3}u^2\Delta y
 +\frac{b_{03}}{2b_{02}}u^3
 +\sum_{i+j\ge4}v_{ij}u^i\Delta y^j.
 \end{aligned}
 \label{vchange}
\end{equation}
Substituting \eqref{vchange} into the first equation of \eqref{po2}, we obtain
\begin{equation}
 \begin{aligned}
 \Delta x={}&\frac{b_{12}}{4b_{02}^2}\Delta y^2
 +\sum_{j\ge3}\hat x_{0j}\Delta y^j
 +u\left[-\frac{b_{12}^2}{2b_{02}^3}\Delta y^2
 +\sum_{j\ge3}\hat x_{1j}\Delta y^j\right]\\
 &+u^2\left[-b_{02}-\frac{3b_{03}}{2b_{02}}\Delta y
 +\sum_{j\ge2}\hat x_{2j}\Delta y^j\right]\\
 &+u^3\left[-\frac23b_{12}
 +\left(-\frac{b_{13}}{b_{02}}
 +\frac{2b_{12}b_{03}}{b_{02}^2}\right)\Delta y
 +\sum_{j\ge2}\hat x_{3j}\Delta y^j\right]\\
 &+\sum_{i\ge4}u^i\sum_{j\ge0}\hat x_{ij}\Delta y^j.
 \end{aligned}
 \label{xfull}
\end{equation}
Hadamard's lemma (see, e.g., \cite[Chap.~II, \S\,6, Lemma~6.10]{golubitsky}) and the inverse function theorem yield an analytic change of variable with the expansion
\begin{equation}
 u=U-\frac{b_{12}}{3b_{02}}U^2
 +\frac{5b_{12}^2}{18b_{02}^2}U^3+\sum_{i\ge4}u_iU^i,
 \label{uchange}
\end{equation}
under which $\Delta x=-b_{02}U^2$ when $\Delta y=0$.

After \eqref{uchange}, relation \eqref{xfull} becomes $\Delta x=G(U,\Delta y)$, where $G$ is analytic.
Under the substitution $U=W+h(\Delta y)$, Taylor expansion in $W$ at fixed $\Delta y$ gives
\[
 \Delta x=G(h(\Delta y),\Delta y)
 +G_U(h(\Delta y),\Delta y)W+O(W^2).
\]
To eliminate the term linear in $W$, define $h(\Delta y)$ by
$G_U(h(\Delta y),\Delta y)=0$.
Since $G_U(0,0)=0$ and $G_{UU}(0,0)=-2b_{02}\ne0$,
the implicit function theorem gives a unique analytic function $h(\Delta y)$ with $h(0)=0$.
Thus the following substitution eliminates the term linear in $W$:
\begin{equation}
 U=W+h(\Delta y).
 \label{hchange}
\end{equation}

The expansion of $h(\Delta y)$ begins with
\[
  h(\Delta y)=-\frac{b_{12}^2}{4b_{02}^4}\Delta y^2
 +\frac{3b_{12}(3b_{12}b_{03}-b_{02}b_{13})}{8b_{02}^6}\Delta y^3
 +O(\Delta y^4).
\]
After the change of variable \eqref{hchange}, we have
\begin{equation}
 \begin{aligned}
 \Delta x={}&\frac{b_{12}}{4b_{02}^2}\Delta y^2
 +\sum_{j\ge3}w_{0j}\Delta y^j\\
 &+W^2\left[-b_{02}-\frac{3b_{03}}{2b_{02}}\Delta y
 +\sum_{j\ge2}w_{2j}\Delta y^j\right]
 +\sum_{i\ge3}W^i\sum_{j\ge1}w_{ij}\Delta y^j.
 \end{aligned}
 \label{xs1}
\end{equation}

Set $s=-\Delta x/b_{02}$. Define the analytic function $s(\Delta y,W)$ by multiplying the right-hand side of \eqref{xs1} by $-1/b_{02}$.
Since $s(0,W)=W^2$, this is a one-parameter deformation of the germ $W^2$
at the origin. The unfolding $R^2+S$ is $R$-versal
(see, e.g., \cite[\S\,8, Corollary~8.4]{arnold1972}). Consequently, the
original relation \eqref{xs1} can be reduced to $R^2+S=0$ by the following
invertible analytic changes of variables:
\begin{equation}
 \begin{gathered}
 R=W\left[1+\sum_{j\ge1}r_{1j}(\Delta y)^j\right]
 +\sum_{i\ge2,\,j\ge1}r_{ij}W^i(\Delta y)^j,\qquad R(0,W)=W,\\
 S=-s+\sum_{i\ge1}s_i(\Delta y)^i,\qquad S(0,0)=0.
 \end{gathered}
 \label{rschange}
\end{equation}
The coefficients $r_{ij}$ and $s_i$ are uniquely determined by the identity
\[
 s(\Delta y,W)=R(\Delta y,W)^2+\sum_{i\ge1}s_i(\Delta y)^i.
\]
The first few coefficients in the transformations \eqref{rschange} are
\[
 s_1=0,\qquad s_2=-\frac{b_{12}}{4b_{02}^3},\qquad
 s_3=\frac{3b_{12}b_{03}-b_{02}b_{13}}{8b_{02}^5},
\]
\[
 r_{11}=\frac{3b_{03}}{4b_{02}^2},\quad
 r_{12}=\frac{3b_{04}}{4b_{02}^3}
 -\frac{27b_{03}^2}{32b_{02}^4}
 -\frac{11b_{12}^3}{24b_{02}^5},\quad
 r_{21}=\frac{b_{02}b_{13}-3b_{12}b_{03}}{2b_{02}^3}.
\]
We thus obtain the canonical fold equation
(a singularity of type $A_2$ \cite[\S\,4, Lemma~4.3]{arnold1972}):
\begin{equation}
 R^2+S=0.
 \label{xs2}
\end{equation}
\begin{samepage}

\begin{theorem}
Let $B=B(u,v)$ be an analytic solution of \eqref{Tre} with the expansion
\eqref{bsol} near $(0,v_*)$, and suppose that $j(0,v_*)=0$ and
$B_{vv}(0,v_*)\ne0$. Then the analytic changes of variables
\eqref{vchange}, \eqref{uchange}, \eqref{hchange}, and \eqref{rschange}
express solutions of the quasilinear system \eqref{tls}, which is equivalent
to the nonlinear stationary transonic equation \eqref{te}, implicitly in
terms of solutions of the canonical fold equation \eqref{xs2}.
\label{foldtheorem}
\end{theorem}
\end{samepage}

\begin{example}
The Tricomi equation \eqref{Tre} has the exact solution
\[
 B=b_{00}+x_*u+y_*(v-v_*)
 +b_{02}\left((v-v_*)^2-\frac{u^3}{3}\right),
\]
where $b_{00},x_*,y_*,b_{02}$ are arbitrary real constants and
$b_{02}\ne0$.
By \eqref{po1}, the changes of variables
\[
 \xi=-\frac{x-x_*}{b_{02}},\qquad
 \eta=v_*+\frac{y-y_*}{2b_{02}}
\]
\nopagebreak[4]
give the fold mapping
\[
 \xi=u^2,\qquad \eta=v.
\]

\end{example}

\section{Genericity of the fold singularity on a compact set}\label{sec-genericity}

This section follows the approach used in \cite[Sec.~4]{sa}, where genericity on compact sets is established for singularities of solutions to quasilinear systems arising in one-dimensional gas dynamics and nonlinear geometrical optics.

Consider analytic mappings $f:(u,v)\mapsto(x,y)$ on an open neighborhood
$U$ of $(0,v_*)$. Equip their space with the Whitney $C^\infty$ topology
(see, e.g., \cite[Chap.~II, \S\,3, Definition~3.1]{golubitsky}).
The first jet extension $j^1f:U\to J^1(U,\mathbb R^2)$ is
\[
 j^1f(u,v)=(u,v,x,y,x_u,x_v,y_u,y_v).
\]
Here $x$, $y$, and their first derivatives are evaluated at $(u,v)\in U$.
The first-jet space is therefore eight-dimensional. Differentiating $j^1f$
involves second derivatives of $x(u,v)$ and $y(u,v)$.

\begin{samepage}
\begin{lemma}
Let $L$ be the closed stratified subset of the first-jet space defined by
\[
 L=\{-u y_v^2-y_u^2=0\}.
\]
The analytic mappings $f$ whose first jet extensions $j^1f$ are transverse
to every stratum of $L$ at points $(0,v)\in U$ form an open dense subset
of the space of analytic mappings. The transversality condition is
\begin{equation}
 y_u(0,v)=0\quad\Longrightarrow\quad y_v(0,v)\ne0.
 \label{translemma}
\end{equation}
\end{lemma}
\end{samepage}
\begin{proof}
The only strata that intersect $u=0$ are the following:
\[
 \begin{array}{ll}
 L_1^\pm=\{-u y_v^2-y_u^2=0,\ \pm y_v>0\},
   &\operatorname{codim}L_1^\pm=1,\\
 L_0=\{u=0,\ y_u=y_v=0\},
   &\operatorname{codim}L_0=3.
 \end{array}
\]
The remaining strata $L_2^\pm=\{\pm u>0,\ y_u=y_v=0\}$ do not meet the line
$u=0$. Whitney's condition~(a) holds for each pair $(L_1^\pm,L_0)$
(see, e.g., \cite[Terminology, condition~$(a)$]{trotman}): the tangent
space to $L_0$ is contained in the tangent spaces to $L_1^\pm$ and in
their limits.
The definition of transversality and Thom's theorem require the strata
to be smooth submanifolds. We verify that the strata $L_1^\pm$ are smooth near their points with $u=0$, and that $j^1f$ is transverse to these strata at such points. It suffices to check that the first partial derivative with
respect to $u$ of the Jacobian expression is nonzero at
$u=0$, $y_u=0$, and $y_v\ne0$:
\[
 \left.\frac{\partial}{\partial u}
 \bigl(-u y_v(u,v)^2-y_u(u,v)^2\bigr)\right|_{u=0}
 =-y_v(0,v)^2\ne0.
\]
The tangent space to $L_0$ has dimension five, whereas the image of
the derivative of $j^1f$ has dimension at most two. Their sum therefore
has dimension at most $5+2=7<8$, so transversality is impossible at a
point $(u,v)\in U$ such that $j^1f(u,v)\in L_0$. If $j^1f(U)$ does not intersect $L_0$, there is no point at which transversality to $L_0$ can fail. Consequently,
$j^1f$ is transverse to $L_0$ if and only if
$j^1f(U)\cap L_0=\varnothing$.

We first prove openness. Condition \eqref{translemma} is equivalent to
$j^1f(U)\cap L_0=\varnothing$. Since $L_0$ is closed, the definition of
the Whitney $C^1$ topology shows that the requirement
$j^1g(u,v)\notin L_0$ for all $(u,v)\in U$ defines an open set of
mappings (see, e.g., \cite[Chap.~II, \S\,3, Definition~3.1]{golubitsky}).
This set is also open in the Whitney $C^\infty$ topology. Its intersection with the space of analytic mappings is open in the induced topology \cite[Chap.~2, \S\,4]{hirsch}.

We now prove density. We show that every open neighborhood of any
analytic mapping $f$ contains a mapping satisfying \eqref{translemma}.
Thom's jet transversality theorem
\cite[Chap.~3, \S\,2, Theorem~2.9]{hirsch} states that the smooth
mappings $g$ for which $j^1g$ is transverse to $L_0$ form a residual
set. Since the space of smooth mappings with the Whitney topology is a Baire space \cite[Chap.~2, \S\,4]{hirsch}, this residual set is dense.
By the definition of the subspace topology, the given neighborhood is
the intersection of an open subset of $C^\infty(U,\mathbb R^2)$
with the subspace of analytic mappings.
Choose within that open set a smooth mapping $f_1$ such that
$j^1f_1(U)\cap L_0=\varnothing$.

Choose an open neighborhood of $f_1(u,v)=(x_1(u,v),y_1(u,v))$ inside the previously chosen open subset of $C^\infty(U,\mathbb R^2)$, such that every mapping in this neighborhood satisfies \eqref{translemma}. The openness proved above guarantees that this choice is possible.
As $U$ is an open subset of $\mathbb R^2$, Grauert's Proposition~8
\cite{grauert} applies separately to the smooth real-valued functions
$x_1(u,v)$ and $y_1(u,v)$. Approximate them by analytic functions
$x_2(u,v)$ and $y_2(u,v)$ closely enough, including their derivatives,
that $f_2=(x_2,y_2)$ belongs to the chosen open neighborhood of $f_1$.
Then $f_2$ belongs to the original neighborhood of $f$ and satisfies
\eqref{translemma}. This proves density.
\end{proof}

For mappings \eqref{po1} associated with solutions of \eqref{Tre},
the expression $-u y_v^2-y_u^2$ equals the Jacobian \eqref{j}.
When $y=B_v$, condition \eqref{translemma} becomes
$B_{uv}(0,v)=0\Rightarrow B_{vv}(0,v)\ne0$.

Choose any compact neighborhood $K\subset U$ of $(0,v_*)$ and set
$I=K\cap\{u=0\}$. We consider the restrictions to $K$ of solutions of
\eqref{Tre} analytic on $U$. Equip this solution space with the Whitney
$C^\infty$ topology on $K$, which agrees with the topology of uniform
convergence on $K$ of derivatives of every order.

\begin{theorem}
In the space of analytic solutions of \eqref{Tre} considered on $K$,
those solutions satisfying
\begin{equation}
 B_{uv}(0,v)=0\quad\Longrightarrow\quad B_{vv}(0,v)\ne0.
 \label{good}
\end{equation}
at every $(0,v)\in I$ form an open dense subset. For each such solution,
the first jet extension of $f=(B_u,B_v)$ is transverse to all strata of
$L$ considered in the lemma at points of $I$. At each point
$(0,v)\in I$ where $j(0,v)=0$, the mapping $f$ has the fold singularity.
\end{theorem}

\begin{proof}
On $u=0$, \eqref{Tre} and \eqref{j} give
\[
 B_{uu}(0,v)=0,\qquad j(0,v)=-B_{uv}(0,v)^2.
\]
At each $(0,v)\in I$ where $j(0,v)=0$, we have
\[
 j_u(0,v)=-B_{vv}(0,v)^2.
\]
By \eqref{good}, $B_{vv}(0,v)\ne0$ at every such point; Theorem~\ref{foldtheorem}
then gives the fold singularity. The lemma also shows that
$j^1f(I)$ does not meet $L_0$ and that $j^1f$ is transverse to
$L_1^\pm$ at points of $I$.

To prove density, we follow the scheme of the proof of the parametric
transversality theorem \cite[Chap.~3, \S\,2, Theorem~2.7]{hirsch}.

Fix an arbitrary analytic solution $B$. Let $\lambda_1$ and $\lambda_2$
be constant real parameters, and set
\[
 B_{\lambda_1,\lambda_2}
 =B+\lambda_1uv+\lambda_2\left(\frac{v^2}{2}-\frac{u^3}{6}\right).
\]
Since both added polynomials satisfy \eqref{Tre}, linearity implies that $B_{\lambda_1,\lambda_2}$ is also an analytic solution, with
\[
 (B_{\lambda_1,\lambda_2})_{uv}(0,v)=B_{uv}(0,v)+\lambda_1,\qquad
 (B_{\lambda_1,\lambda_2})_{vv}(0,v)=B_{vv}(0,v)+\lambda_2.
\]
Failure of condition \eqref{good} requires the following two equalities to hold at the same point $(0,v)\in I$:
\[
 \lambda_1=-B_{uv}(0,v),\qquad \lambda_2=-B_{vv}(0,v).
\]
Every such parameter pair belongs to the image of the smooth mapping
\begin{equation}
 v\longmapsto\bigl(-B_{uv}(0,v),-B_{vv}(0,v)\bigr).
 \label{badparameters}
\end{equation}
Its domain is the open subset $\{v:(0,v)\in U\}$ of the real line, and
its codomain is the $(\lambda_1,\lambda_2)$-parameter plane. Its derivative
matrix is
\[
 \begin{pmatrix}-B_{uvv}(0,v)\\-B_{vvv}(0,v)\end{pmatrix}.
\]
The derivative has rank at most one, whereas the parameter space is two-dimensional. Every point in the image of this mapping is therefore a critical value. The mapping is smooth because $B$ is analytic. Sard's theorem \cite[Chap.~3, \S\,1, Theorem~1.3]{hirsch} therefore implies that its image has two-dimensional measure zero. Every open ball centered at $(0,0)$
has positive measure, so it contains a parameter pair outside this
image. Such a pair can be chosen arbitrarily close to $(0,0)$ and
satisfies \eqref{good} throughout $I$.

It remains to show that these perturbed solutions can be made arbitrarily close to $B$ in the chosen topology. For any integers $p,q\ge0$,
\[
 \sup_K\bigl|\partial_u^p\partial_v^q(B_{\lambda_1,\lambda_2}-B)\bigr|
 \longrightarrow0\quad\text{as }(\lambda_1,\lambda_2)\to(0,0).
\]
Indeed, the difference of the solutions depends linearly on the
parameters, and all derivatives of the added polynomials are bounded
on the compact set $K$. Given any neighborhood of $B$ in the solution space, all parameter pairs in a sufficiently small open ball centered at $(0,0)$ produce solutions in that neighborhood.

Choose a pair in the ball but outside the image of \eqref{badparameters},
which has measure zero. For this pair, condition \eqref{good} holds
throughout $I$: if both $(B_{\lambda_1,\lambda_2})_{uv}$ and $(B_{\lambda_1,\lambda_2})_{vv}$ vanished at some point of $I$, the parameter pair would belong to that image.

Thus every neighborhood of every solution $B$ contains a solution
satisfying \eqref{good}. This proves density.

To prove openness, observe that \eqref{good} is equivalent to
$B_{uv}(0,v)^2+B_{vv}(0,v)^2>0$ at every point of the compact set $I$.
Consequently,
\[
 m=\min_{(0,v)\in I}\bigl(B_{uv}(0,v)^2+B_{vv}(0,v)^2\bigr)>0.
\]
Continuity of the derivatives and compactness of $I$ also give
\[
 M=\max_{(0,v)\in I}\bigl(|B_{uv}(0,v)|+|B_{vv}(0,v)|\bigr)<\infty.
\]
Choose $0<\delta<\min\{1,m/[4(M+1)]\}$ and consider solutions
$\widetilde B$ satisfying
\[
 \sup_K|\widetilde B_{uv}-B_{uv}|<\delta,\qquad
 \sup_K|\widetilde B_{vv}-B_{vv}|<\delta.
\]
These inequalities define an open neighborhood of $B$ in the
solution space under consideration. In what follows, all derivatives
are evaluated at $(0,v)\in I$. Write
$\widetilde B_{uv}=B_{uv}+(\widetilde B_{uv}-B_{uv})$ and square both
sides to obtain
\[
 \widetilde B_{uv}^2-B_{uv}^2
 =2B_{uv}(\widetilde B_{uv}-B_{uv})
 +(\widetilde B_{uv}-B_{uv})^2.
\]
As $I\subset K$, we have
\[
 |\widetilde B_{uv}^2-B_{uv}^2|
 \le2|B_{uv}|\delta+\delta^2.
\]
Likewise, $|\widetilde B_{vv}^2-B_{vv}^2|\le2|B_{vv}|\delta+\delta^2$.
The triangle inequality now yields
\[
 \left|\widetilde B_{uv}^2+\widetilde B_{vv}^2-B_{uv}^2-B_{vv}^2\right|
 \le2M\delta+2\delta^2<2(M+1)\delta<\frac m2.
\]
Since $B_{uv}^2+B_{vv}^2\ge m$ on $I$, the preceding estimate gives $\widetilde B_{uv}^2+\widetilde B_{vv}^2>m/2>0$.
The two derivatives cannot vanish simultaneously, so \eqref{good}
holds for every solution in the constructed neighborhood.
\end{proof}

\section{Concluding remarks}\label{sec-remarks}

\begin{remark}
If the compact set $K$ is replaced by a bounded open rectangle $U$
containing $(0,v_*)$, the solutions satisfying \eqref{good} on
$U\cap\{u=0\}$ are not dense in the strong Whitney $C^\infty$ topology.
Consider the following basic neighborhood of the zero solution: analytic solutions $B$ of \eqref{Tre} whose values and derivatives up to order two have absolute values less than $\operatorname{dist}(z,\partial U)$ at each $z\in U$
\cite[Chap.~II, \S\,3, pp.~42--43]{golubitsky}.

Multiplying \eqref{Tre} by $B_u$ and integrating by parts yields
$\tfrac12\iint_U B_v^2\,du\,dv=0$: the boundary terms vanish because
$B_u$ and $B_v$ tend to zero as $(u,v)$ approaches $\partial U$. Thus $B$ is independent
of $v$. Since $B$ itself also tends to zero at the boundary,
$B\equiv0$. Hence this neighborhood contains only the zero solution,
which fails to satisfy \eqref{good}.

The same obstruction to extending the genericity arguments from a
compact set to an open neighborhood applies to the systems of
one-dimensional gas dynamics and nonlinear geometrical optics
\cite[Theorems~1--3 and Proposition~1]{sa}. Let $V$ be an open
rectangle, and equip the space of solutions on $V$ with the strong
Whitney topology.

Consider a neighborhood of the zero solution in which the absolute
values of $B$ and its derivatives up to order two are less than the
distance to the boundary of $V$. If $V$ adjoins the axis $\rho=0$,
the distance is measured to the three outer sides; no decay to zero
is required at the axis itself. For the hyperbolic equation
$B_{rl}=\kappa(B_r-B_l)$ (gas dynamics, $\rho_*>0$), choose a small
rectangle on which $\kappa$ is nonzero. Integrating the identity
$\partial_l(B_r^2)-\partial_r(B_l^2)=2\kappa(B_r-B_l)^2$
gives $B_r=B_l=0$. For Chaplygin gas, $\kappa=0$ and $B_{rl}=0$,
so $B=f(r)+g(l)$. At fixed $r$, the derivative $B_r=f'(r)$ is
independent of $l$. Letting $l$ approach a side of the rectangle
gives $f'(r)=0$; similarly, $g'(l)=0$.

For the elliptic equation
$\rho B_{\rho\rho}+2B_\rho=-\alpha(\rho)B_{uu}$
(nonlinear geometrical optics, $\rho_*>0$), where $\alpha>0$,
integration by parts gives
$\iint_V(\rho^2B_\rho^2+\rho\alpha B_u^2)\,du\,d\rho=0$.
When $\rho_*=0$, the equation is
$\rho B_{\rho\rho}+2B_\rho=\pm\alpha(\rho)B_{uu}$.
Multiplying it by $\rho^3B_\rho$ and integrating over a small
rectangle in the half-plane $\rho\ge0$, we obtain
$\iint_V\rho^2(3\alpha+\rho\alpha')B_u^2\,du\,d\rho=0$;
here $3\alpha+\rho\alpha'>0$.

In all four cases, these identities and the boundary behavior of $B$
imply $B\equiv0$. The zero solution fails the conditions of the
corresponding theorems and Proposition~1. Consequently, their
density results on compact sets do not extend to these open rectangles
in the strong Whitney topology.

This remark was written with the assistance of ChatGPT~6.2.
\end{remark}

\begin{remark}
In the proof of Theorem~2, there is no need to rule out
$\lambda_1=\lambda_2=0$. For this parameter pair,
$B_{\lambda_1,\lambda_2}=B$. If $B$ already satisfies \eqref{good}
on $I$, the zero pair is admissible.

Otherwise, $B_{uv}(0,v)$ and $B_{vv}(0,v)$ vanish simultaneously
at some $(0,v)\in I$. The pair $(0,0)$ then belongs to the image
of \eqref{badparameters}. Choosing parameters outside the image of \eqref{badparameters} therefore automatically excludes $(0,0)$.
\end{remark}

\begin{remark}
In the elliptic region ($u>0$), by \eqref{j}, the condition $j=0$
implies $B_{uv}=B_{vv}=0$. Equation \eqref{Tre} then also gives
$B_{uu}=0$ at the same point, so the mapping \eqref{po1} has corank two
there. This motivates the conjecture that one can describe only a
two-dimensional section of the elliptic umbilic ($D_4^{-}$) in the
elliptic region, whereas in the hyperbolic region ($u<0$) one can
describe the fold ($A_2$), cusp ($A_3$), and a two-dimensional section
of the hyperbolic umbilic ($D_4^{+}$).
\end{remark}

\begin{remark}
In the approach adopted here
\cite{ks1995,suleimanov2009,sul12,chap24,ss2020,ss22,ss24,mss24,shav,mss2026,sa},
we use only relations among Taylor coefficients that follow from
\eqref{Tre} and from the vanishing of the Jacobian \eqref{j}.
Further equality constraints on the coefficients in \eqref{bsol}
are not imposed; they lie outside the scope of this approach.
Expand the solution $B$ of \eqref{Tre} at $(u_*,v_*)$.
In terms of its Taylor coefficients $b_{ij}$, the condition
$j(u_*,v_*)=0$ reads
\[
 -4u_*b_{02}^2-b_{11}^2=0.
\]
If $u_*>0$, then $b_{02}=b_{11}=0$; if $u_*<0$, then
$b_{11}=\pm2\sqrt{-u_*}\,b_{02}$.
For $u_*=0$, we necessarily have $b_{11}=0$, while
$b_{20}=0$ follows from \eqref{Tre}.

For the solution $B=B(u,v)$ of \eqref{Tre}, the coefficients $b_{ij}$ are determined by derivatives evaluated at $(u_*,v_*)$. Thus the displayed
equalities constrain the location of $(u_*,v_*)$ in the $(u,v)$-plane.
With $u_*=0$ fixed, the condition $b_{11}=0$ becomes a condition on $v_*$.
The relation $b_{02}=0$ follows neither from \eqref{Tre} nor
from the vanishing of \eqref{j}.

If one goes beyond this approach and imposes the additional conditions
$b_{02}=0$ and $b_{12}\ne0$, a two-dimensional section of a $D_5$
singularity can be described. Indeed, in the generating family
\[
 F=-B(u,v)+xu+yv
\]
(whose critical points with respect to $u$ and $v$ are precisely the solutions
of \eqref{po1}, and which is introduced by analogy with
\cite[Introduction]{suleimanov2009}), make the invertible analytic
changes of variables
\[
 \begin{aligned}
 u={}&R-\frac{b_{03}}{b_{12}}L
 -\frac{13b_{03}^{\,2}}{6b_{12}^{\,2}}R^2
 +\left(\frac{7b_{03}^{\,3}}{3b_{12}^{\,3}}
 -\frac{b_{13}}{b_{12}}\right)RL\\
 &+\left(-\frac{5b_{03}^{\,4}}{6b_{12}^{\,4}}
 +\frac{b_{03}b_{13}}{b_{12}^{\,2}}
 -\frac{b_{04}}{b_{12}}\right)L^2,\\
 \Delta v={}&L+\frac{b_{03}}{6b_{12}}R^2.
 \end{aligned}
\]
The Jacobian determinant of this change of variables at the origin is $1$.
When $\Delta x=\Delta y=0$,
\[
 F=-b_{00}+y_*v_*-b_{12}RL^2+\frac{b_{12}}6R^4
 +O\bigl((|R|+|L|)^5\bigr).
\]
Since $b_{12}\ne0$ and the 4-jet of $F$ is sufficient
\cite[\S\,3, Remark~3.5; \S\,5, Lemma~5.4]{arnold1972},
the germ of $F$ at $\Delta x=\Delta y=0$ has a singularity of type
$D_5$. A two-dimensional section can be chosen if the derivative
matrix of two parameters of a versal unfolding with respect to
$(\Delta x,\Delta y)$ is nonsingular at the origin. The coefficients
of the terms linear in $R,L$ form the pair
$\bigl(\Delta x,\ \Delta y-(b_{03}/b_{12})\Delta x\bigr)$,
whose Jacobian determinant with respect to $(\Delta x,\Delta y)$ is $1$.
The parameters multiplying nonlinear terms in the normal form can
therefore be expressed in terms of those multiplying the linear terms.
\end{remark}

\begin{remark}
The result also applies to other physically relevant models whose governing systems either have the form \eqref{tls} or can be reduced to it. For
example, the equations governing longitudinal motion of an elastic
rod \cite[formula~(1.1)]{rohde2005}, with the quadratic stress law
$\sigma(a)=a^2/2$, are
\[
 \begin{cases}
  a_\tau-q_s=0,\\
  q_\tau-aa_s=0,
 \end{cases}
\]
Here $a=a(s,\tau)$ is the longitudinal strain, $q=q(s,\tau)$ the
velocity, $s$ the material coordinate, and $\tau$ time.
The substitutions $x=s$, $y=\tau$, $u=-a$, $v=-q$ reduce this system
to \eqref{tls}.

Another example is the asymptotic gas-dynamical system near the
caustics of weak shock waves \cite[Sec.~IV, formula~(63)]{rosalestabak1998}:
\[
 \begin{cases}
  \bigl(YA-A^2/2\bigr)_X+Q_Y=0,\\
  A_Y+Q_X=0.
 \end{cases}
\]
Here $A=A(X,Y)$ and $Q=Q(X,Y)$ are normalized perturbations of
density and transverse velocity, and $X,Y$ are scaled coordinates
near the caustic. The substitutions $x=X$, $y=Y$, $u=A-Y$,
$v=-Q-X$ reduce this system to \eqref{tls}.
\end{remark}

\begin{remark}
One possible physical interpretation of the fold singularity described
here is that, as a gradient-catastrophe point is approached, the
magnitude of the acceleration of fluid particles grows without bound
while their velocity remains finite.
In \cite[Sec.~Criterion for the Breakdown of Isentropic Irrotational Flow]{tsien1944},
unbounded acceleration and an unbounded pressure gradient are described as
characteristic features of limiting lines in a flow. The possibility
that isentropic irrotational flow breaks down before such a line is
reached is also discussed. In \cite[Introduction]{tsienkuo1946},
the appearance of a limiting line was associated with the breakdown
of isentropic irrotational flow. For the families of flows around
profiles of bounded curvature considered in \cite[Chap.~5, \S\,20]{bers2016},
however, the emergence of shock waves cannot be explained by the
sudden appearance of a limiting line. The monograph
\cite[Chap.~20, Sec.~20.5]{shapiro1954} discusses two possible outcomes:
the formation of shock waves or the absence of a stationary flow.
The possible connection between the fold singularity described in this paper and these physical phenomena requires further study.
\end{remark}

\section*{Acknowledgements}
The author thanks B.~I.~Suleimanov for helpful discussions.

\section*{Statements and Declarations}

\subsection*{Funding}
This research was supported by the Russian Science Foundation,
grant No.~26-11-00021,
\url{https://rscf.ru/project/26-11-00021/}.

\subsection*{Competing interests}
The author declares that there are no competing interests.

\subsection*{Data availability}
No datasets were generated or analyzed during this study.

\subsection*{Author contributions}
The author contributed to all aspects of this work.

\end{document}